\documentclass[12pt]{amsart}

\usepackage{amssymb,amsmath,setspace,mathrsfs}
\newtheorem{proposition}{Proposition}[section]
\newtheorem{lemma}[proposition]{Lemma}
\newtheorem{theorem}[proposition]{Theorem}

\theoremstyle{remark}
\numberwithin{equation}{section}
\def\a{\alpha}
\def\de{\delta}
\def\ga{\gamma}

\def\o{\omega}

\def\dimp{{\rm dim}_{_{\rm P}}}
\def\Dim{{\rm Dim}}
\def\uB{\hbox{${\rm B}$-$\overline{\dim}$}}

\def\Pim{\hbox{${\rm P}$-$\dim$}}
\def\DimB{\hbox{${\rm B}$-$\overline{\rm dim}$}}

\def\I{\mathbf{1}}
\def\R{{\mathbf R}}
\def\N{{\mathbf R}^N}
\def\d{{\mathbf R}^d}

\def\E{{\mathrm E}}
\def\P{{\mathrm P}}

\def\fBM{{$f\text{BM}$}}

\begin{document}
\begin{spacing}{1.35}

\title[Dimension-Profiles and  \fBM\ ]{%
    Packing-Dimension Profiles\\ and Fractional Brownian Motion}
\thanks{Research partially supported by NSF grant DMS-0404729.}
\author[D. Khoshnevisan]{Davar Khoshnevisan}
\address{Department of Mathematics, 155 S. 1400 E., JWB 233,
        University of  Utah,
        Salt Lake City, UT 84112--0090}
\email{davar@math.utah.edu}
\urladdr{http://www.math.utah.edu/\~{}davar}

\author{Yimin Xiao}
\address{Department of Statistics and Probability, A-413 Wells
        Hall, Michigan State University,
        East Lansing, MI 48824}
\email{xiao@stt.msu.edu}
\urladdr{http://www.stt.msu.edu/\~{}xiaoyimi}

\date{November 10, 2006}

\maketitle

\begin{abstract}
    In order to compute the packing dimension of orthogonal projections
    Falconer and Howroyd (1997) introduced a family of
    packing dimension profiles $\Dim_s$ that are parametrized by
    real numbers $s>0$. Subsequently,
    Howroyd (2001) introduced alternate $s$-dimensional
    packing dimension profiles $\Pim_s$ and proved,
    among many other things, that $\Pim_s E=\Dim_s E$
    for all integers $s>0$ and all analytic sets $E\subseteq\R^N$.

    The goal of this article is to prove that
    $\Pim_s E=\Dim_s E$ for all real numbers
    $s>0$ and analytic sets $E\subseteq\R^N$.
    This answers a question of Howroyd (2001, p.\ 159).
    Our proof hinges on a new property of fractional Brownian motion.
\end{abstract}

\section{Introduction}

Packing dimension and packing measure were introduced in the early
1980s by Tricot (1982) and Taylor and Tricot (1985) as dual concepts to
Hausdorff dimension and Hausdorff measure. Falconer
(1990) and Mattila (1995) contain systematic accounts.

It has been known for some time now that some Hausdorff dimension
formulas --- such as those for orthogonal projections and those
for image sets of fractional Brownian motion
--- do not have packing dimension analogues; see J\"arvenp\"a\"a (1994)
and Talagrand and Xiao (1996) for precise statements. This suggests
that a new concept of dimension is needed to compute
the packing dimension of some random sets.

In order to compute the packing dimension of orthogonal
projections Falconer and Howroyd (1997) introduced a family of
\emph{packing dimension profiles} $\{\Dim_s\}_{s>0}$ that we
recall in Section \ref{sec:profile} below. Falconer and Howroyd
(1997) proved that for every analytic set $E \subset \R^N$ and
every integer $1 \le m \le N$,
\begin{equation}\label{Eq:FH97}
    \dimp \left( P_{_V} E \right)
    = \Dim_{m} E \qquad \hbox{ for $\gamma_{n, m}$-almost
    all $V \in \mathscr{G}_{n, m}$},
\end{equation}
where $\gamma_{n, m}$ is the natural orthogonally-invariant measure
on the Grassman manifold $\mathscr{G}_{n, m}$ of all $m$-dimensional
subspaces of $\R^N$, and $P_{_V}E$ denotes the projection of $E$
onto $V$.

Subsequently, Howroyd (2001) introduced a family $\{\uB_s\}_{s>0}$
of box-dimension profiles, together with their regularizations
$\{\Pim_s\}_{s>0}$. The latter are also called packing dimension
profiles; see Section \ref{sec:profile}. Howroyd (2001) then used
these dimension profiles to characterize the [traditional] box and
packing dimensions of orthogonal projections. In addition, Howroyd
(2001, Corollary 32) proved that for all analytic sets
$E\subseteq\R^N$: (i) $\Pim_s E \ge \Dim_{s} E$ if $s > 0$; and
(ii) if $s \in (0\,, N)$ is an \emph{integer} then
\begin{equation} \label{Eq:Howroyd1}
    \Pim_s E = \Dim_{s} E.
\end{equation}
Finally, $\Pim_s E$ and $\Dim_{s} E$ agree for arbitrary
$s \ge N$, and their common value is the packing dimension $\dimp E$.

The principle aim of this note is to prove that \eqref{Eq:Howroyd1}
holds for all real numbers $s\in(0\,,N)$. Equivalently,
we offer the following.

\begin{theorem}\label{th:main}
    Equation \eqref{Eq:Howroyd1} is valid for all $s>0$.
\end{theorem}
This solves a question of Howroyd (2001, p.\ 159).

Our derivation is probabilistic, and relies on properties
of fractional Brownian motion (\fBM). In order to explain
the connection to  \fBM\  let $X :=\{X(t)\}_{t \in \R^N}$ be a
$d$-dimensional  \fBM\  with Hurst
parameter $H\in(0\,,1)$. That is,
$X(t) =( X_1(t)\,, \ldots, X_d(t))$ for all $t\in\R^N$,
where $X_1, \ldots, X_d$ are independent copies of a real-valued
 \fBM\  with common Hurst parameter $H$
(Kahane, 1985, Chapter 18).
Xiao (1997) proved that for every analytic set $E \subseteq \N$,
\begin{equation} \label{Xiao97}
    \dimp X(E) = \frac{1}{H}\, \Dim_{Hd} E \quad \ \hbox{ a.s.}
\end{equation}
Here we will derive an alternative expression.

\begin{theorem}\label{th:fBM}
    For all analytic sets $E \subseteq \R^N$,
    \begin{equation} \label{KX1}
        \dimp X(E) = \frac1 {H}\,
        \Pim_{Hd} E \quad  \ \hbox{ a.s.}
    \end{equation}
\end{theorem}

Thanks to (\ref{Xiao97}) and Theorem \ref{th:fBM},
$\Dim_{Hd}E=\Pim_{Hd}E$ for all integers $d\ge 1$
and all $H\in(0\,,1)$. Whence follows Theorem \ref{th:main}.

We establish Theorem \ref{th:fBM} in
Section 3, following the introductory Section 2 wherein we introduce
some of stated notions of fractal geometry in
greater detail. Also we add a Section 4 where
we derive yet another equivalent formulation for
the $s$-dimensional packing dimension profile $\Dim_s E$
of an analytic set $E\subseteq\N$. We hope to
use this formulation of $\Dim_s E$ elsewhere
in order to compute the
packing dimension of many interesting random sets.

Throughout we will use the letter $K$ to denote an unspecified
positive and finite constant whose value may
differ from line to line and sometimes even within the same line.

\section{Dimension Profiles}
\label{sec:profile}
In this section  we recall briefly aspects of the
theories of dimension profiles of Falconer and Howroyd (1997) and
Howroyd (2001).
\subsection{Packing Dimension via Entropy Numbers}
For all $r > 0$ and all bounded sets $E \subseteq
{\R^N}$ let $N_r (E)$ denote the maximum number of
disjoint closed balls of radius  $r$ whose respective
centers are all in $E$. The [upper]
\emph{box dimension} of $E$ is defined as
\begin{equation}\label{Eq:DimB}
    \DimB\, E = \limsup_{r \downarrow 0} \frac{\log N_r (E)} {\log
    (1/r)}.
\end{equation}

We follow Tricot (1982) and define
the \emph{packing dimension} of $E$ as
the ``regularization'' of $\DimB\, E$. That is,
\begin{equation}\label{Eq:Tricot}
    \dimp E = \inf \left\{\sup_{k\ge 1}\, \DimB\, F_k:  \ E \subseteq
    \bigcup_{k=1}^\infty F_k\right\}.
\end{equation}

There is also a corresponding notion of the packing dimension
of a Borel measure. Indeed, the [lower] \emph{packing dimension}
of a Borel measure $\mu$ on $\R^N$ is
\begin{equation}\label{Eq:MuDim}
    \dimp \mu = \inf \left\{ \dimp E:\
    \mu (E) > 0 \ \hbox{and $E \subseteq
    \R^N$ is a Borel set}\right\}.
\end{equation}
One can compute $\dimp E$  from $\dimp \mu$ as well:
Given an analytic set $E \subseteq \N$ let
$\mathscr{M}^+_c(E)$ denote the collection of
all finite compactly-supported Borel measures on $E$. Then,
according to Hu and Taylor (1994),
\begin{equation} \label{Eq:dimpset}
    \dimp E = \sup \left\{ \dimp \mu :\ \mu \in \mathscr{M}^+_c(E)
    \right\}.
\end{equation}

\subsection{The Packing Dimension Profiles of Falconer and Howroyd}
Given a finite Borel measure $\mu$ on $\N$ and an $s\in(0\,,\infty]$ define
\begin{equation}
    F_s^{\mu} (x\,, r) := \int_{\N} \psi_s
    \left(\frac{x-y}{r}\right)
    \, \mu (dy),
\end{equation}
where for finite $s\in(0\,,\infty)$,
\begin{equation}\label{Eq:psi}
    \psi_s (x) := \min \left( 1 \,, |x|^{-s} \right) \qquad
    {}^\forall \, x \in \R^N,
\end{equation}
and $\psi_\infty := \mathbf{1}_{\{y\in\R^d:\, |y|\le 1\}}$.
The \emph{$s$-dimensional packing dimension
profile} of $\mu$ is defined as
\begin{equation}\label{Def:Dimprof}
    \Dim_{s}\mu = \sup \left\{ t \ge 0:\
    \liminf_{r \downarrow 0}\
    \frac{F_s^{\mu} (x\,, r)}{r^t}
    = 0\  \hbox{ for  $\mu$-a.a.\ $x \in \R^N$}
    \right\}.
\end{equation}
Packing dimension profiles generalize
the packing dimension because
$\dimp \mu = \Dim_{s}\mu$ for all
finite Borel measures $\mu$ on $\R^N$ and for all $s\ge N$.
See Falconer and Howroyd (1997, p. 272) for a proof.

Falconer and Howroyd (1997) also
defined the $s$-dimensional
packing dimension profile of a Borel set $E
\subseteq \N$ by
\begin{equation} \label{Eq:dimprofileset}
    \Dim_s E = \sup\left\{ \Dim_s \mu :\
    \mu \in \mathscr{M}^+_c(E) \right\}.
\end{equation}

\subsection{The Packing Dimension Profiles of Howroyd}
If $E \subset \R^N$ and $s > 0$, then a sequence
of triples $(w_i\,, x_i\,, r_i)_{i= 1}^\infty$
is called a {\it $(\psi_s\,, \de)$-packing of
$E$} whenever $w_i \ge 0$, $x_i \in E$, $0 < r_i \le \de$, and
\begin{equation}
    \sup_{i\ge 1}
    \sum_{j=1}^\infty w_j\, \psi_s\bigg(\frac {x_i - x_j} {r_j}\bigg)
    \le 1.
\end{equation}

For all $E \subset \R^N$, define
\begin{equation}\label{Def:KerPre}
    \mathcal{P}^{\a, s}_0(E) := \lim_{\de\downarrow 0}\,
    \sup \left\{ \sum_{i=1}^\infty w_i\, (2r_i)^\a\, :\,
    (w_i\,, x_i\,,
    r_i)_{i=1}^\infty  \ \hbox{ is a $(\psi_s\,, \de)$-packing of $E$}
    \right\}.
\end{equation}
Then the \emph{$\a$-dimensional $\psi_s$-packing measure
$\mathcal{P}^{\a, s}(E)$} is defined as
\begin{equation}\label{Def:kerpack}
    \mathcal{P}^{\a, s}(E) = \inf \left\{ \sum_{k=1}^\infty
    \mathcal{P}^{\a, s}_0(E_k):\, E \subseteq \bigcup_{k=1}^\infty E_k
    \right\}.
\end{equation}
The \emph{$s$-dimensional
packing dimension profile of $E$} can then be defined as
\begin{equation}\label{Def:Pim}
    \Pim_s E := \inf\left\{ \alpha > 0:\, \mathcal{P}^{\a, s}(E)=
    0\right\}.
\end{equation}

We will make use of the following two lemmas.
They are ready consequences of
Lemma 20 and Theorem 22 of Howroyd (2001), respectively.

\begin{lemma}\label{Lem:nonsigma1}
    If $E \subset \R^N$ and  $\mathcal{P}^{\ga, s}(E)> 0$, then
    $E$ has non-sigma-finite
    $\mathcal{P}^{\a, s}$-measure  for
    every $\a \in (0\,, \ga)$.
\end{lemma}

\begin{lemma}\label{Lem:nonsigma}
    Let $A \subset \R^N$ be an analytic set of non-sigma-finite
    $\mathcal{P}^{\a, s}$-measure. Then there exists a compact set $K
    \subset A$ such that $\mathcal{P}^{\a, s}_0 (K \cap G) = \infty$
    for all open sets $G \subset \R^N$ with $K \cap G \ne \varnothing$.
    Moreover, $K$ is also of non-sigma-finite
    $\mathcal{P}^{\a,s}$-measure.
\end{lemma}

\subsection{Upper Box Dimension Profiles}
Given $r > 0$ and $E \subset \R^N$, a sequence of pairs
$(w_i\,, x_i)_{ i=1}^k$ is a \emph{size-$r$ weighted $\psi_s$-packing
of $E$} if: (i) $x_i \in E$; (ii) $w_i \ge 0$; and (iii)
\begin{equation}\label{Eq:psipacking}
    \max_{1\le i\le k}
    \sum_{j=1}^k w_j \psi_s \left( \frac{x_i - x_j}{r}
    \right) \le 1.
\end{equation}
Define
\begin{equation}\label{Eq:Npsi}
    N_r(E\,; \psi_s) := \sup
    \left\{ \sum_{i=1}^k w_i:\ (w_i\,, x_i)_{i=1}^k
    \hbox{ is a  size-$r$ weighted $\psi_s$-packing of $E$}
    \right\}.
\end{equation}
This quantity is related to the entropy number $N_r(E)$. In fact,
Howroyd (2001, Lemma 5) has shown  that $N_r(E\,; \psi_\infty) =
N_{r/2}(E)$ for all $ r> 0$ and all $E\subseteq \R^N$. We will use
this fact in the proof of Lemma \ref{Lem:Meas} below.

The \emph{$s$-dimensional upper box dimension of $E$} is defined as
\begin{equation}\label{Def:uBpsi}
    \uB_s \, E := \limsup_{r \downarrow 0}
    \frac{\log N_r(E\,; \psi_s)} {\log(1/r)},
\end{equation}
where $\log 0 := - \infty$. Note in particular that $\uB_s \, \varnothing = -
\infty.$  It is possible to deduce that $s \mapsto \uB_s \, E$ is
non-decreasing.

Define $\mathscr{P}_A(E)$ to be the collection of all probability
measures that are supported on a finite number
of points in $E$. For all $\mu\in
\mathscr{P}_A(E)$ define
\begin{equation}\label{Eq:J}
    J_s \left(r\,, \mu \right)
    := \max_{x\in\text{supp}\,\mu}\,F_s^{\mu} (x\,, r)
    \quad\text{and}\quad
    I_s\left( r\,, \mu \right)
    := \int F_s^{\mu} (x\,, r)  \, \mu(dx).
\end{equation}
For $E \subset \R^N$, define
\begin{equation}\label{Eq:Z}
    Z_s(r\,; E) := \inf_{\mu\in\mathscr{P}_A(E)} J_s \left( r\,, \mu\right).
\end{equation}
Howroyd (2001) has demonstrated that for all $s,r>0$,
\begin{equation}\label{Eq:Z2}
    Z_s(r\,; E) = \inf_{\mu\in\mathscr{P}_A(E)}
        I_s \left( r\,, \mu\right)\quad\text{and}\quad
    N_r \left( E\,; \psi_s \right) = \frac{1}{Z_s(r\,; E)}.
\end{equation}
Consequently,
\begin{equation}\label{Def:uBpsi2}
    \uB_s \, E = \limsup_{r \downarrow 0}
    \frac{\log Z_s(r\,; E)} {\log r}.
\end{equation}

According to Howroyd (2001, Proposition 8),
\begin{equation} \label{Eq:Boxdim}
    \uB_s \, E = \DimB \, E\qquad
    {}^\forall s\ge N,\, E\subseteq\R^N.
\end{equation}
Howroyd (2001) also proved that $\Pim_s$ is
the regularization of $\uB_s$; i.e.,
\begin{equation}\label{Eq:Box-Packing}
    \Pim_s E = \inf \left\{\sup_{k\ge 1} \, \uB_s\, E_k: \ E \subseteq
    \bigcup_{k=1}^\infty E_k \right\},
\end{equation}
This is the dimension-profile analogue of \eqref{Eq:Tricot}.

\section{Proof of Theorem \ref{th:fBM}}
\label{sec:fBM}

Recall that $X$ is a centered, $d$-dimensional, $N$-parameter Gaussian
random field such that for all $s,t\in\R^N$
and $j,k\in\{1\,,\ldots,d\}$,
\begin{equation} \label{Eq:Ycov}
    \text{Cov}\left( X_j (s) \,, X_k(t) \right)=
    \frac12
    \left(|s|^{2H}+|t|^{2H}-|s-t|^{2H}\right)\delta_{ij}.
\end{equation}
Throughout, we assume that the process $X$ is constructed in a
complete probability space $(\Omega\,, {\mathscr F}, \P)$, and
that $t\mapsto X(t\,,\omega)$ is continuous
for almost every $\omega \in \Omega$.
According to the general theory of Gaussian processes
this can always be arranged.

Our proof of Theorem \ref{th:fBM} hinges on several lemmas. The
first is a technical lemma which verifies the folklore statement
that, for every $r>0$ and $E \subseteq \R^N$, the entropy number
$N_r(X(E))$ is a random variable. We recall that
$(\Omega\,,\mathscr{F},\P)$ is assumed to be complete.
\begin{lemma}\label{Lem:Meas}
    Let $E \subseteq \R^N$ be a fixed set,
    and choose and fix some $r > 0$. Then  $N_r(X(E))$ and
    $Z_\infty(r\,; X(E)) $ are non-negative random variables.
\end{lemma}

\begin{proof}
    It follows from (\ref{Eq:Z2})  that $Z_\infty(r\,; X(E)) =
    1/N_{r/2}\left(X(E)\right)$. Hence it suffices to prove
    $N_r(X(E))$ is a random variable.

    Let $C(\R^N)$ be the space of continuous functions $f:
    \R^N \to \R^d$ equipped with the norm
    \begin{equation}
        \|f\| = \sum_{k=1}^\infty 2^{-k} \frac{\max_{|t|\le k} |f(t)|} { 1
        + \max_{|t| \le k} |f(t)|}.
    \end{equation}
    According to general theory we can assume
    without loss of generality that $\Omega=C(\R^N)$.
    It suffices to prove that for all $a>0$ fixed,
    $\Theta_a := \{f \in C(\R^N): N_r(f(E))>
    a\}$ is open  and hence Borel measurable.
    For then $\{\omega\in\Omega:\ N_r(X(E))>a\}
    = X^{-1}(\Theta_a)$
    is also measurable.

    To this end we assume that $N_r(f(E))> a$, and define $n := \lfloor a
    \rfloor$.  There necessarily exist $t_1, \ldots,
    t_{n+1}\in E$ such that
    $|f(t_i) - f(t_j) | > 2r$ for all
    $1 \le i \ne j \le n+1$.
    Choose and fix $\eta \in (0\,, 1)$ such that $\eta <
    \min\{|f(t_i) - f(t_j) | - 2r: \ {}^\forall\, 1 \le i \ne j \le
    n+1\}$. We can then find an integer $k_0>0$ such that
    $|t_i | \le k_0$ for
    all $ i = 1, \ldots, n+1$. It follows from our
    definition of the norm $\|\cdot\|$ that for all $\delta \in
    \big( 0\,, \eta\,2^{-(k_0+2)} \big)$ and all functions $g \in C(\R^N)$
    with $\|g-f\|< \delta$,
    \begin{equation}
        \max_{1\le i\le n+1} \left| g(t_i) - f(t_i)
        \right| < \frac{\eta}{2}.
    \end{equation}
    This and the triangle inequality
    imply $|g(t_i) - g(t_j)| \ge |f(t_i) - f(t_j)| - \eta > 2r$
    for all $1 \le i \ne j \le n+1$, and hence
    $N_r(g(E)) > n$. This verifies that $\{f
    \in C(\R^N): N_r(f(E))> a\}$ is an open set.
\end{proof}

The following lemma is inspired by Lemma 12 of Howroyd (2001).
We emphasize that
$\E \left[ Z_\infty(r\,; X(E))
\right]$ is well defined (Lemma \ref{Lem:Meas}).

\begin{lemma}\label{Lem:LB1}
    If  $E\subseteq\R^N$ then
    \begin{equation}
        \E[ Z_\infty(r\,; X(E))] \le
        K\, Z_{Hd} \left( r^{1/H}\!; E \right)
        \qquad{}^\forall r > 0.
    \end{equation}
    The constant $K \in(0\,,\infty)$ depends only on $d$ and $H$.
\end{lemma}

\begin{proof}
    Note that $(\mu\circ X^{-1})\in\mathscr{P}_A(X(E))$
    whenever  $\mu\in\mathscr{P}_A(E)$.
    Hence, $Z_\infty(r\,; X(E))\le I_\infty (r\,,\mu\circ X^{-1})$.
    Because $I_\infty( r\,, \mu \circ X^{-1})
    = \iint \I_{\{|X(s) - X(t)| \le r\}}\, \mu(ds)\, \mu(dt)$
    for all $r>0$,
    \begin{equation}\label{Eq:I2}\begin{split}
        \E \left[ Z_\infty(r\,; X(E))\right]
            &\le \iint \P\left\{|X(s) - X(t)| \le r\right\}
            \,\mu(ds)\, \mu(dt)\\
        &\le K \iint \left( \frac {r^d} {|s-t|^{H d}}
            \wedge 1\right)\, \mu(ds)\, \mu(dt)
            = K\, I_{Hd}
            \left(r^{1/H},\, \mu\right),
    \end{split}\end{equation}
    where the last inequality follows from the self-similarity and stationarity of
    the increments of $X$, and where $K>0$ is a constant
    that depends only on $d$ and $H$.
    We obtain the desired result by
    optimizing over all $\mu\in\mathscr{P}_A(E)$.
\end{proof}

\begin{lemma}\label{Lem:LB2}
    For all nonrandom sets $E \subset \R^N$,
    \begin{equation}\label{Eq:Box}
        \DimB \, X(E) \ge \frac 1 H\, \uB_{Hd}\, E \quad \hbox{ a.s.}
    \end{equation}
\end{lemma}

\begin{proof}
    Without loss of generality we assume $\uB_{Hd} E > 0$,
    for otherwise there is nothing left to prove.
    Then for any constant $\gamma \in (0\,, \uB_{Hd} E)$ there
    exists a sequence $\{r_n\}_{n=1}^\infty$ of positive numbers such that $r_n
    \downarrow 0$ and $Z_{Hd} (r_n\,; E) = o(r_n^\gamma)$
    as $n\to\infty$.
    It follows from Lemma \ref{Lem:LB1} and
    Fatou's lemma  that
    \begin{equation}
        \E \left[ \liminf_{r \downarrow 0}
        \frac{ Z_\infty \left(r\,;  X(E) \right)}{r^{\gamma/H}} \right]
        \le \liminf_{n \to \infty} \frac{\E\left[Z_\infty \left(r_n^H\,;
        X(E) \right)\right]}{r_n^{\gamma}}
        \le K\, \lim_{n \to \infty}
        \frac{Z_{Hd} (r_n\,; E)}{r_n^{\gamma}} = 0.
    \end{equation}
    Consequently, (\ref{Def:uBpsi2}) and (\ref{Eq:Boxdim})
    together imply that $\DimB\, X(E) \ge \gamma/H$ a.s.
    The lemma follows because
    $\gamma \in (0\,,  \uB_{Hd} E)$ is arbitrary.
\end{proof}

The following Lemma is borrowed from  Falconer and
Howroyd (1996, Lemma 5).

\begin{lemma}\label{Lem:Boxpacking}
    If a set $E \subset \R^N$ has the property that $\DimB (E \cap G)
    \ge \de$ for all open sets $G \subset \R^N$ such that $E \cap G \ne
    \varnothing$. Then $\dimp E \ge \de$.
\end{lemma}

We are ready to prove Theorem \ref{th:fBM}.

\begin{proof}[Proof of Theorem \ref{th:fBM}]
    Since $\Pim_{H d} E
    \ge \Dim_{H d} E$, \eqref{Xiao97} implies that $ \dimp X(E)$
    is almost surely bounded above by
    $\frac1 {H}\, \Pim_{H d} E$. Consequently, it remains to
    prove the reverse inequality.

    To this end we may assume without loss of generality
    that $\Pim_{H d} E> 0$, lest the inequality becomes
    vacuous. Choose and fix an arbitrary $\alpha
    \in (0\,,\Pim_{H d} E)$. Lemma \ref{Lem:nonsigma1} implies that
    $E$ has non-$\sigma$-infinite $\mathcal{P}^{\alpha, Hd}$-measure. By
    Lemma \ref{Lem:nonsigma}, there exists a compact set $K \subset E$
    such that $\mathcal{P}^{\alpha, Hd} (G\cap K) = \infty$ for all open
    sets $G \subset \R^N$ with $G\cap K \ne \varnothing$.

    By separability there exists a
    countable basis of the usual euclidean topology on $\R^N$. Let $\{G_k\}_{%
    k = 1}^\infty$ be an enumeration of those sets in the basis that
    intersect $K$. It follows from Lemma \ref{Lem:LB2} that for every
    $k=1, 2,\ldots$ there exists an event $\Omega_k$ of
    $\P$-measure one
    such that for all $\omega \in \Omega_k$,
    \begin{equation}\label{Eq:Ok}
        \DimB \, X_\o(G_k \cap K) \ge \frac 1 H \uB_{Hd} (G_k \cap K) \ge
        \frac \a H.
    \end{equation}
    Therefore, $\Omega_0 := \bigcap\limits_{k=1}^\infty\Omega_k$
    has full $\P$-measure, and for every $\omega \in \Omega_0$,
    \begin{equation}\label{Eq:Ok2}
        \DimB \, \big(X_\o(K)\cap U \big)\ge \DimB \, \big(X_\o(K\cap
        X^{-1}(U) \big) \ge \frac \a H.
    \end{equation}
    The preceding is valid for
    all open sets $U$ with $X(K)\cap U \ne \varnothing$ because
    $X^{-1}(U) $ is open and $K\cap X^{-1}(U) \ne \varnothing$.
    According to Lemma \ref{Lem:Boxpacking}  this proves that
    $\dimp X_\o(K) \ge \a/H$ almost surely.
    Because $\a \in (0\,, \Pim_{Hd} E)$ is arbitrary
    this finishes the proof of
    Theorem \ref{th:fBM}.
\end{proof}

\section{An Equivalent Definition}

Given a Borel set $E\subset\R^N$, we define $\mathscr{P}(E)$ as the
collection all probability measures $\mu$ on $\R^N$ such that  $\mu(E) = 1$
[$\mu$ is called a probability measure on $E$].
Define for all Borel sets $E\subset\R^N$ and all $s\in (0\,,\infty]$,
\begin{equation}\label{Eq:Z3}
    \mathscr{Z}_s (r\,;E) := \inf_{\mu\in\mathscr{P}(E)}
    I_s( r\,,\mu).
\end{equation}
Thus, the sole difference between $\mathscr{Z}_s$ and $Z_s$ is
that in the latter we use all finitely-supported [discrete]
probability measures on $E$, whereas in the former we use all
probability measures on $E$. We may also define  $\mathscr{Z}_s$
using $\mathscr{M}^+_c(E)$ in place of $\mathscr{P}(E)$ in
(\ref{Eq:Z3}). Our next theorem shows that all these notions lead
to the same $s$-dimensional box dimension.

\begin{theorem}\label{th:reform}
    For all analytic sets $E\subset\R^N$
    and all $s\in (0\,,\infty]$,
    \begin{equation}
    \uB_s(E) = \limsup_{r\downarrow 0}
    \frac{\log \mathscr{Z}_s(r\,;E)}{\log r}.
    \end{equation}
\end{theorem}

\begin{proof}
    Because $\mathscr{P}_A(E)\subset\mathscr{P}(E)$ it follows
    immediately that $\mathscr{Z}_s(r\,;E) \le Z_s(r\,;E)$.
    Consequently,
    \begin{equation}
        \limsup_{r\downarrow 0}
        \frac{\log \mathscr{Z}_s(r\,;E)}{\log r}
        \ge \uB_s(E).
    \end{equation}

    We explain the rest only when $N=1$; the general
    case is handled similarly.

    Without loss of much generality suppose $E\subset [0,1)$ and $\mu$ is
    a probability measure on $E$.

    For all integers $n\ge 1$ and $i\in\{0\,,1\,,\ldots,n-1\}$ define
    $C_i=C_{i,n}$ to be $1/n$ times the half-open interval
    $[ i\,,i+1)$.
    Then, we can write
    $I_s ( 1/n \,, \mu ) = T_1 + T_2$,
    where
    \begin{equation}\begin{split}
        T_1 & := \mathop{\sum\sum}_{\substack{0\le i<n \\
            j\in\{i-1,i,i+1\}}}
            \int_{C_i}\int_{C_j} \left( 1
            \wedge \frac{1}{n|x-y|} \right)^s\, \mu(dx)\, \mu(dy),\\
        T_2& := \mathop{\sum\sum}_{\substack{0\le i<n \\
            j\not\in\{i-1,i,i+1\}}}
            \int_{C_i}\int_{C_j} \left( 1
            \wedge \frac{1}{n\,|x-y|} \right)^s\, \mu(dx)\, \mu(dy).
    \end{split}\end{equation}

    Any interval $C_j$ with $\mu(C_j)=0$ does not contribute
    to $I_s ( 1/n \,, \mu )$. For every $j$ with $\mu(C_j) > 0$, we choose
    an arbitrary  point $\tau_j \in E\cap C_j$  and denote $w_j := \mu(C_j)$.
    Then the discrete probability measure $\nu$ that puts mass
    $w_j$ at $ \tau_j\in E$ belongs to $\mathscr{P}_A(E)$. For simplicity
    of notation, in the following we assume $\mu(C_j) > 0$ for all $j=0\,,1\,,\ldots,n-1$.

    If $j\not\in\{i-1,i,i+1\}$, then
    $\sup_{x\in C_i} \sup_{y\in C_j} |x-y| \le  3\,|\tau_j- \tau_i|$,
    whence we have
    \begin{equation}
        T_2 \ge \frac{1}{3^s}\mathop{\sum\sum}_{\substack{0\le i<n \\
        j\not\in\{i-1,i,i+1\}}}
        \left( 1 \wedge \frac{1}{n\, |\tau_j-\tau_i|} \right)^s
        w_i\, w_j.
    \end{equation}
    If $j\in\{i-1\,,i\,,i+1\}$, then a similar case-by-case analysis
    can be used. This leads us to the bound,
    \begin{equation}\label{eq:BBBB}\begin{split}
        I_s\left( \frac{1}{n} ~,~ \mu\right)
            &\ge \frac{1}{3^s}\mathop{\sum\sum}_{0\le i,j<n}
            \left( 1 \wedge \frac{1} {n\,|\tau_j-\tau_i|} \right)^s
            w_i\, w_j\\
        &= \frac{1}{3^s}\iint \left( 1\wedge \frac{1/n}{|a-b|} \right)^s
            \, \nu(da)\, \nu(db).
    \end{split}\end{equation}
    Consequently, the right-hand side of \eqref{eq:BBBB}
    is at most $3^{-s}\,Z_s(1/n\,; E)$.
    It follows that
    \begin{equation}
    3^{-s}\, Z_s\left (\frac 1 n \,; E \right) \le \mathscr{Z}_s \left( \frac 1 n \,; E \right) \le
    Z_s \left( \frac 1 n \,; E \right).
    \end{equation}
    If $r$ is between $1/n$ and $1/(n+1)$, then
    $Z_s(r\,;E)$ is between $Z_s(1/n\,;E)$
    and $Z_s(1/(n+1)\,;E)$. A similar remark applies
    to $\mathscr{Z}_s$. Because $\log n\sim
    \log(n+1)$ as $n\to\infty$,
    this proves the theorem.
\end{proof}

\bigskip

\bibliographystyle{plain}
\begin{small}

\end{small}
\end{spacing}
\end{document}